\documentclass[a4paper,11pt,leqno]{amsart}

\RequirePackage[T1]{fontenc}
\RequirePackage[utf8]{inputenc}
\RequirePackage{xcolor}
	\definecolor{lapis}{HTML}{26619C}
\RequirePackage{amsfonts,amssymb,amsmath}
\RequirePackage{amsthm}
\RequirePackage{mathtools}
\RequirePackage{tikz-cd}
\RequirePackage{soul}

\RequirePackage{hyperref}
\hypersetup{colorlinks,allcolors=lapis}
\RequirePackage{lmodern,stmaryrd,microtype}

\usepackage[capitalize]{cleveref}

\numberwithin{equation}{section}
\theoremstyle{plain}
	\newtheorem{thm}{Theorem}[section]
	\newtheorem{conj}[thm]{Conjecture}
	
\theoremstyle{definition}
	\newtheorem{ex}[thm]{Example}
\theoremstyle{remark}	
	\newtheorem{rem}[thm]{Remark}

\Crefname{thm}{Theorem}{Theorems}
\Crefname{conj}{Conjecture}{Conjectures}
\Crefname{ex}{Example}{Examples}

\usepackage{autonum}

\newcommand{\bb}{\mathbf}
\newcommand{\ur}{\mathrm}

\newcommand{\erase}[1]{\textcolor{white}{#1}}
\newcommand{\lapis}[1]{\textcolor{lapis}{#1}}
	\newcommand{\dfemph}[1]{\lapis{\emph{#1}}}

\DeclareMathOperator{\End}{End}

\newcommand{\Br}{\mathit{Br}}
\newcommand{\Deg}{\mathrm{Deg}}
\newcommand{\GL}{\mathrm{GL}}
\newcommand{\Irr}{\mathrm{Irr}}
\newcommand{\Uch}{\mathrm{Uch}}

\newcommand{\core}[1]{\textit{$#1$-cor}}

\newcommand{\rat}{\mathrm{rat}}

\calclayout

\makeatletter
\g@addto@macro \normalsize {%
	\setlength\abovedisplayskip{10pt plus 2pt minus 2pt}%
	\setlength\belowdisplayskip{10pt plus 2pt minus 2pt}%
}
\makeatother
\title{Level-Rank Dualities for Finite Reductive Groups}

\author[Minh-T\^{a}m Trinh]{Minh-T\^{a}m Quang Trinh}
\address{Department of Mathematics, Howard University, Washington, DC 20059}
\email{minhtam.trinh@howard.edu}

\author{Ting Xue}
\address{School of Mathematics and Statistics, University of Melbourne, VIC 3010, Australia}
\email{ting.xue@unimelb.edu.au}
\thanks{TX was supported in part by the ARC grant DP250100824.}

\begin{document}

\begin{abstract}
This is an extended abstract of our work ``Level-Rank Dualities from $\Phi$-Cuspidal Pairs\ldots''
We present evidence for a family of surprising coincidences within the representation theory of a finite reductive group $G$: more precisely, dualities between blocks of cyclotomic Hecke algebras attached by Brou\'e--Malle to $\Phi$-cuspidal pairs of $G$, where the Hecke parameters are specialized not to the order of the underlying finite field, but to roots of unity.
For the groups $G = \mathrm{GL}_n(\mathbf{F}_q)$, these coincidences can be expressed very concretely in terms of the combinatorics of partitions, and the whole story recovers an avatar of the level-rank duality studied by Frenkel, Uglov, Chuang--Miyachi, and others.
\end{abstract}

\maketitle

\thispagestyle{empty}



\setcounter{section}{-1}

\section{Introduction} 

This extended abstract of our work~\cite{tx} emphasizes motivations from combinatorics, modular representation theory, and classical level-rank duality, while omitting any detailed discussion of affine Springer fibers.

As we will explain in \S\ref{subsec:level-rank}, the \emph{level-rank duality} discovered by Igor Frenkel in the early 1980s was a surprising reciprocal relationship in the representation theory of affine Lie algebras, between a level-$l$ integrable module over $\widehat{\mathfrak{sl}}_m$ and a level-$m$ integrable module over $\widehat{\mathfrak{sl}}_l$.
By work of Ariki, this relationship, or rather, its deformation quantization, can be recast in terms of the modular representations of Hecke algebras attached to the wreath products $S_N \ltimes \bb{Z}_m^N$, where $\bb{Z}_m$ is $\bb{Z}/m\bb{Z}$ and $S_N$ is the symmetric group on $N$ letters.
In this form, it can be upgraded to an equivalence of derived categories by work of Rouquier--Shan--Varagnolo--Vasserot.

A major motivation for the wreath products above is that they are relative Weyl groups for certain twisted Levi subgroups of the finite general linear groups $\GL_n(\bb{F}_q)$.
The main contribution of~\cite{tx} is to recast level-rank duality itself, so that it can be generalized to a family of conjectures about an arbitrary finite reductive group $G$.
The key ingredient is a generalization of Harish--Chandra theory due to Brou\'e--Malle--Michel, building on work of Lusztig, which furnishes the Hecke algebras that generalize Ariki's.
Our strongest conjecture is a derived equivalence relating the modular representations of a \emph{$\Phi_m$-cuspidal} Hecke algebra at an $l$th root of unity and a \emph{$\Phi_l$-cuspidal} Hecke algebra at an $m$th root of unity.
Here, the prefixes $\Phi_l$- and $\Phi_m$- roughly encode the extent to which the relevant Levi subgroups of $G$ are twisted: \emph{i.e.}, fail to be split.

We hope that our work contributes a new richness to the work of Brou\'e--Malle--Michel on finite reductive groups, and to the modular representation theory of these groups in general.
Following the appearance of the first version of~\cite{tx}, Chlouveraki--Malle extended some of our results for exceptional groups, and proposed spetsial generalizations of some of our conjectures~\cite{cm_zeitschrift}.

Although our strongest conjecture involves categories, we stress that our other conjectures lead to concrete numerical predictions, which one can check by hand or computer.
Moreover, in the proof that our categorical conjecture recovers the theorem of Rouquier--Shan--Varagnolo--Vasserot, the key tool is a purely elementary combinatorial structure: a bijection $\Upsilon_m^l$ between sets of charged multipartitions, due to Uglov.

We start by reviewing the combinatorics of $\Upsilon_m^l$ in \Cref{sec:partitions}.
In \Cref{sec:reps}, we review Ariki's Hecke algebras and their modular representations, then explain level-rank duality in these terms.
In \Cref{sec:cats}, we explain how the duality is categorified.
In \Cref{sec:reductive}, we turn to finite reductive groups $G$, reviewing the theories of Lusztig and Brou\'e--Malle--Michel; in \Cref{sec:linear}, we show how they connect to the concepts from \Crefrange{sec:partitions}{sec:reps} when $G = \GL_n(\bb{F}_q)$.
Finally, we state our conjectures and results, illustrated with examples, in \Cref{sec:tx}.

We thank the referee for useful comments and corrections to the text.
We declare that AI was not used in any part of the preparation of this article.

\section{Partitions}\label{sec:partitions}

\subsection{}

Let $\Pi$ be the set of all integer partitions.
We will regard a partition as a weakly decreasing sequence $\pi = (\pi_1, \pi_2, \ldots)$ such that $\pi_i = 0$ for all $i$ large enough; in practice, we write a partition like $\pi = (3, 3, 1, 0, \ldots)$ using a shorthand like $(3^2, 1)$.
We define the size of $\pi$ to be $|\pi| \vcentcolon= \pi_1 + \pi_2 + \ldots$ and the length of $\pi$ to be the number of nonzero entries.

For any positive integer $m$, we define an \dfemph{$m$-partition} to be an $m$-tuple of partitions $\vec{\pi} = (\pi^{(0)}, \pi^{(1)}, \ldots, \pi^{(m - 1)})$, indexed from $0$ through $m - 1$.
We will occasionally regard these indices as elements of $\bb{Z}_m \vcentcolon= \bb{Z}/m\bb{Z}$.
We define the size of $\vec{\pi}$ to be $|\vec{\pi}| \vcentcolon= |\pi^{(0)}| + \cdots + |\pi^{(m - 1)}|$.

\subsection{}

The combinatorial part of our story starts with an analogue, for partitions, of the notion of long division by $m$.
Recall that a partition is an \dfemph{$m$-core} if and only if it contains no hook lengths divisible by $m$.
For instance, the only $1$-core is the empty partition; the $2$-cores are the ``staircase'' partitions $1$, $(2, 1)$, $(3, 2, 1)$, \emph{etc.}; but $m$-cores for $m \geq 3$ are more complicated.
Let $\Pi_{\core{m}} \subseteq \Pi$ be the subset of $m$-cores.
Then the analogue of the map sending an integer to its quotient and remainder modulo $m$ is a bijection $\Pi \xrightarrow{\sim} \Pi^m \times \Pi_{\core{m}}$, usually known as the \dfemph{core-quotient bijection}.

Following Uglov~\cite{uglov}, we define it thus:
Let $\bb{B}$ be the collection of subsets $\beta \subseteq \bb{Z}$ with the property that $n \in \beta$ when $n$ is a sufficiently negative integer, but $n \notin \beta$ when $n$ is sufficiently positive.
There is a bijection $\Pi \times \bb{Z} \xrightarrow{\sim} \bb{B}$ that takes a pair $(\pi, s)$ to the set
\begin{align}
	\beta_{\pi, s} = \{\pi_i - i + s \mid i = 1, 2, \ldots\}.
\end{align}
There is also a bijection $\upsilon_m = (\upsilon_m^{(0)}, \upsilon_m^{(1)}, \ldots, \upsilon_m^{(m - 1)}) \colon \bb{B} \xrightarrow{\sim} \bb{B}^m$, defined by setting 
\begin{align}
	\upsilon_m^{(r)}(\beta) = \{q \in \bb{Z} \mid mq + r \in \beta\}
	\quad\text{for $0 \leq r < m$}.
\end{align}
The composition
\begin{align}
	\Upsilon_m \colon \Pi \times \bb{Z} \xrightarrow{\sim} \bb{B} \xrightarrow{\upsilon_m} \bb{B}^m \xrightarrow{\sim} (\Pi \times \bb{Z})^m = \Pi^m \times \bb{Z}^m
\end{align}
recovers the core-quotient map as follows:
The $m$-quotient of $\pi$ is the $m$-partition formed by the $\Pi^m$-component of $\Upsilon_m(\pi, 0)$.
The $m$-core of $\pi$ is the unique partition $\mu$ such that $\Upsilon_m(\mu, s)$ and $\Upsilon_m(\pi, s)$ have the same $\bb{Z}^m$-component for any $s$.

Note that $\mu$ is an $m$-core if and only if the $\Pi^m$-component of $\Upsilon_m(\mu, s)$ is the empty $m$-partition for some (or any) $s$.
Also, note that there is a purely graphical way to calculate the core-quotient map, using Young diagrams~\cite[\S{2.7}]{jk}.

In the literature, elements of $\Pi^m \times \bb{Z}^m$ are called \dfemph{charged $m$-partitions} and written with the ket notation $|\vec{\pi}, \vec{s}\rangle$ of quantum mechanics.
The vector $\vec{s}$ is called the \dfemph{$m$-charge}.
Elements of $\bb{B}^m$ are sometimes called \dfemph{$m$-runner abacus configurations}, or \dfemph{$m$-abaci}.
They can be pictured as subsets of $\bb{Z} \times \bb{Z}_m$, or more vividly, as configurations of beads on a horizontal abacus with $m$ runners.
Under the identification $\Pi^m \times \bb{Z}^m \simeq \bb{B}^m$, the operation $|\vec{\pi}, \vec{s}\rangle \mapsto |\vec{\emptyset}, \vec{s}\rangle$ corresponds to sliding all beads as far left as possible.

\begin{ex}\label[ex]{ex:3-abacus}
	Let $\pi = (8, 6, 1)$ and $s = 2$.
	Then $\beta_{\pi, s} \in \bb{B}$ is represented by the following $1$-abacus:
	\begin{align}
		\small
		\begin{array}{@{\quad\cdots\quad}cccccccccccccc}
			-4 &-3 &-2 &\erase{-}- &\erase{-}0 &\erase{-}- &\erase{-}- &\erase{-}- &\erase{-}- &\erase{-}- &\erase{-}6 &\erase{-}- &\erase{-}- &\erase{-}9
		\end{array}
	\end{align}
	So $\Upsilon_3(\pi, s) \in \bb{B}^3$ is represented by the following $3$-abacus, where we have labeled each position $(q, r)$ with the value of $3q + r$.
	Note that $r$ increases \emph{downward}.
	\begin{align}
		\small
		\begin{array}{@{\quad\cdots\quad}cccccccc}
			-12 &-9 &-6 &-3 &\erase{-}0 &\erase{-}- &\erase{-}6 &\erase{-}9\\
			-11 &-8 &-5 &-2 &\erase{-}- &\erase{-}- &\erase{-}- &\\
			-10 &-7 &-4 &\erase{-}- &\erase{-}- &\erase{-}- &\erase{-}- &
		\end{array}
	\end{align}
	Sliding the beads at $6$ and $9$ as far left as possible gives the $3$-abacus of $\Upsilon_3(\mu, s)$, where $\mu = (5, 3, 1)$.
	So the $3$-core of $\pi$ is $(5, 3, 1)$.
\end{ex}

\subsection{}

We will discuss generalizations of $\upsilon_m$ and $\Upsilon_m$ that involve fixing another integer $l > 0$.
Let $\upsilon_m^l = (\upsilon_m^{l, (0)}, \upsilon_m^{l, (1)}, \ldots, \upsilon_m^{l, (m - 1)}) \colon \bb{B}^l \xrightarrow{\sim} \bb{B}^m$ be defined by the rule:
\begin{align}
\begin{array}{r@{\:}c@{\:}ll}
lq + r_l \in \upsilon_m^{l, (r_m)}(\vec{\beta})
&\iff
&mq + r_m \in \beta^{(r_l)}
&\text{for all $q, r_l, r_m \in \bb{Z}$ such that }\\
&&&\text{$0 \leq r_l < l$ and $0 \leq r_m < m$}.
\end{array}
\end{align}
It is also possible to rewrite $\upsilon_m^l$ as a composition $\bb{B}^l \xrightarrow{\upsilon_l^{-1}} \bb{B} \xrightarrow{\upsilon_m^\ast} \bb{B}^m$, where $\upsilon_m^\ast$ is some modified version of $\upsilon_m$.
We can now form the composition
\begin{align}
	\Upsilon_m^l \colon \Pi^l \times \bb{Z}^l = \bb{B}^l \xrightarrow{\upsilon_m^l} \bb{B}^m = \Pi^m \times \bb{Z}^m.
\end{align}
Note that $\upsilon_m^1 = \upsilon_m$ and $\Upsilon_m^1 = \Upsilon_m$.

\begin{ex}\label[ex]{ex:4-abacus}
	The bijection $\Upsilon_4^3$ takes the $4$-abacus
	\begin{align}
		\small
		\begin{array}{@{\quad\cdots\quad}cccccc}
			-12 &-8 &-4 &\erase{-}0 \\
			-11 &-7 &-3 &\erase{-}- \\
			-10 &-6 &-2 &\erase{-}2 \\
			-9 &-5 &\erase{-}- &\erase{-}3
		\end{array}
	\end{align}
	to the $3$-abacus in \Cref{ex:3-abacus}. 
	
	Here is a graphical interpretation of $\Upsilon_4^3$.
	If we subdivide the $4$-abacus above into $4 \times 3$ matrices---grouping together the columns where $q_4 = 0, 1, 2$, the columns where $q_4 = 3, 4, 5$, and so on---and subdivide the $3$-abacus from \Cref{ex:3-abacus} into $3 \times 4$ matrices in an analogous way, then $\Upsilon_4^3$ can be viewed as taking the transpose of each $4 \times 3$ matrix to get a $3 \times 4$ matrix.
	In our examples:
	\begin{align}
		\small
		\begin{array}{ccc@{\quad}|@{\quad}ccc}
			\bullet &\bullet &\bullet &\bullet &- &-\\
			\bullet &\bullet &\bullet &- &- &-\\
			\bullet &\bullet &\bullet &\bullet &- &-\\
			\bullet &\bullet &- &\bullet &- &-
		\end{array}
		\quad\mapsto\quad
		\begin{array}{cccc@{\quad}|@{\quad}cccc}
			\bullet &\bullet &\bullet &\bullet &\bullet &- &\bullet &\bullet \\
			\bullet &\bullet &\bullet &\bullet &- &- &- &-\\
			\bullet &\bullet &\bullet &- &- &- &- &-
		\end{array}
	\end{align}
\end{ex}

In the next section, we will explain how $\Upsilon_m^l$ holds meaning in the representation theory of the wreath products
\begin{equation}\label{eq:wreath}
S_{N, m} \vcentcolon= S_N \ltimes \bb{Z}_m^N,
\end{equation}
or rather, the modular representation theory of associated Hecke algebras.
Note that the $m = 1$ and $m = 2$ cases recover the Weyl groups of types $A$ and $BC$, respectively.

\section{Representations}\label{sec:reps}

\subsection{}

Recall that $S_N$ is a reflection group generated by the system of simple reflections $s_1, s_2, \ldots, s_{N - 1}$, where $s_i$ transposes letters $i$ and $i + 1$.
Writing $\Br_N$ for the braid group on $N$ strands, we have a quotient map $\Br_N \to S_N$, sending the $i$th positive simple twist $\sigma_i$ to the transposition $s_i$.
It descends to a quotient map of rings $H_N(x) \to \bb{Z}S_N$, where the domain is the \dfemph{Hecke algebra}
\begin{align}
	H_N(x) \vcentcolon= \frac{\bb{Z}[x^{\pm 1}]\Br_N}{\langle (\sigma_i - 1)(\sigma_i + x) \mid i = 1, \ldots, N - 1\rangle}.
\end{align}
For any nonzero complex number $\xi$, we set 
\begin{align}
	H_N(\xi) = \bb{C} \otimes_{\bb{Z}[x^{\pm 1}]} H_N(x),
	\quad\text{where $\bb{Z}[x^{\pm 1}] \to \bb{C}$ sends $x \mapsto \xi$}.
\end{align}
For sufficiently generic $\xi  \in \bb{C}^\times$, there is an isomorphism $H_N(\xi) \simeq \bb{C}S_N$, and hence, a bijection between simple $H_N(\xi)$-modules up to isomorphism and irreducible characters of $S_N$.
By work of Frobenius, Schur, and Young, the latter are also in bijection with integer partitions of size $N$.
But $H_N(\xi)$ may not even be semisimple at special $\xi$.

\subsection{}

Now fix a positive integer $m$.
The group $S_{N, m}$ defined by~\eqref{eq:wreath} is a \emph{complex} reflection group, generated by $S_N$ and a pseudo-reflection of order $m$.
The map $\Br_N \to S_N = S_{N, 1}$ generalizes to a map $\Br_{N, m} \to S_{N, m}$, where
\begin{align}
	\Br_{N, m}
	&= \left\{\begin{array}{ll}
		\Br_N
		&m = 1,\\
		\langle \Br_N, \tau \mid
		\text{$\tau\sigma_1\tau\sigma_1 = \sigma_1\tau\sigma_1\tau$ and $\tau \sigma_j = \sigma_j \tau$ for $j \neq 1$}
		\rangle.
		&m > 1.
	\end{array}\right.
\end{align}
The Hecke algebra generalizes to the so-called \dfemph{Ariki--Koike algebra}
\begin{align}
	H_{N, m}(\vec{x})
	= \frac{\bb{Z}[\vec{x}^{\pm 1}]\Br_{N, m}}{
		\left\langle\!\begin{array}{l}
			\text{$(\sigma_i - 1)(\sigma_i + x_\sigma)$ for $i = 1, \ldots, N - 1$},\\
			(\tau - x_{\tau, 0})(\tau - x_{\tau, 1}) \cdots (\tau - x_{\tau, m - 1})
		\end{array}\!\right\rangle},
\end{align}
where $\vec{x}$ denotes a collection of indeterminates $(x_\sigma, x_{\tau, 0}, \ldots, x_{\tau, m - 1})$, and $\vec{x}^{\pm 1}$ means we also adjoin their inverses.
For any $\vec{\xi} \in (\bb{C}^\times)^m$, we define $H_{N, m}(\vec{\xi})$ similarly to how we defined $H_N(\xi)$.

Like before, there is an isomorphism $H_{N, m}(\vec{\xi}) \simeq \bb{C}S_{N, m}$ at generic $\vec{\xi}$, inducing a bijection between simple $H_{N, m}(\vec{\xi})$-modules up to isomorphism and irreducible characters of $S_{N, m}$.
By work of Clifford, the latter are indexed by $m$-partitions of size $N$.
Like before, $H_{N, m}(\vec{\xi})$ can degenerate at special $\vec{\xi}$.

We are especially interested in parameters of the form 
\begin{align}\label{eq:ariki-koike-params}
	\vec{\xi} = (\zeta, \zeta^{\vec{s}}) \vcentcolon= (\zeta, \zeta^{s^{(0)}}, \zeta^{s^{(1)}}, \ldots, \zeta^{s^{(m - 1)}}),
\end{align}
where $\zeta$ is a root of unity and $\vec{s} = (s^{(0)}, s^{(1)}, \ldots, s^{(m - 1)})$ is an $m$-charge.
It turns out that the representation theory of $H_{N, m}(\zeta, \zeta^{\vec{s}})$ encodes and is encoded by the combinatorics of $m$-cores and $m$-quotients.

\subsection{}

Given an associative algebra $H$, we write $\mathsf{K}_0(H)$ for the Grothendieck group of the category of finitely-generated $H$-modules, and $\mathsf{K}_0^+(H) \subseteq \mathsf{K}_0(H)$ for its non-negative part.
Let $H_{N, m}^\circ(\vec{x}) = \bb{K} \otimes_{\bb{Z}[\vec{x}^{\pm 1}]} H_{N, m}(\vec{x})$ for some field $\bb{K} \supseteq \bb{Z}[\vec{x}^{\pm 1}]$ such that $H_{N, m}^\circ(\vec{x}) \simeq \bb{K}S_{N, m}$.

As explained in~\cite[Ch.\@ 7]{gp} and~\cite[Ch.\@ 3]{gj}, the specialization map from $H_{N, m}(\vec{x})$ to $H_{N, m}(\vec{\xi})$ induces a \dfemph{decomposition map}
\begin{align}\label{eq:decomposition-map}
	\mathrm{d}_{\vec{\xi}} \colon \mathsf{K}_0^+(H_{N, m}^\circ(\vec{x})) \to \mathsf{K}_0^+(H_{N, m}(\vec{\xi})).
\end{align}
At the same time, every character of $S_{N, m}$ corresponds to a module over $H_{N, m}^\circ(\vec{x})$ up to isomorphism.
Altogether, we get maps
\begin{align}
	\{\vec{\pi} \in \Pi^m \mid |\vec{\pi}| = N\} \xrightarrow{\sim} \Irr(S_{N, m}) \to \mathsf{K}_0^+(H_{N, m}^\circ(\vec{x})) \xrightarrow{\mathrm{d}_{\vec{\xi}}} \mathsf{K}_0^+(H_{N, m}(\vec{\xi})).
\end{align}
The elements of the image of the composition are the classes of the so-called \dfemph{Specht modules}.
We say that two $m$-partitions of size $N$, or irreducible characters of $S_{N, m}$, belong to the same \dfemph{$\vec{\xi}$-block} if and only if the corresponding Specht modules fit into a finite sequence of $H_{N, m}(\vec{\xi})$-modules in which consecutive modules always share a Jordan--H\"older factor.
Note that $H_{N, m}(\vec{\xi})$ is semisimple if and only if each block is a singleton.
In general, $\vec{\xi}$-blocks can be more complicated.

\subsection{}

Fix positive integers $l, m$.
Fix a primitive $l$th root of unity $\eta_l$ and a primitive $m$th root of unity $\eta_m$.

In a slogan, $\Upsilon_m^l$ relates $(\eta_m, \eta_m^{\vec{r}})$-blocks of $l$-partitions, as we run over $l$-charges $\vec{r}$, and $(\eta_l, \eta_l^{\vec{s}})$-blocks of $m$-partitions, as we run over $m$-charges $\vec{s}$.
To make this precise, we define a \dfemph{rank-$m$, level-$l$ Uglov datum} to be a triple $(K, \vec{r}, \mathbf{a})$, where $K \geq 0$ is any integer, $\vec{r}$ is an $l$-charge, and $\mathbf{a}$ is an $(\eta_m, \eta_m^{\vec{r}})$-block of the $l$-partitions of size $K$.
There is a natural map from the set of charged $l$-partitions to the set of rank-$m$, level-$l$ Uglov data:
To $|\vec{\pi}, \vec{r}\rangle$, we associate the Uglov datum in which $K = |\vec{\pi}|$ and $\mathbf{a}$ is the $(\eta_m, \eta_m^{\vec{r}})$-block containing $\vec{\pi}$.
The following result is implicit in~\cite{uglov}, via a theorem of Ariki~\cite[Thm.\@ 6.2.21]{gj}:

\begin{thm}[Uglov]\label[thm]{thm:uglov}
The bijection $\Upsilon_m^l$ descends to a bijection
\begin{align}
\{\text{rank-$m$, level-$l$ Uglov data}\}
	\xrightarrow{\sim} 
	\{\text{rank-$l$, level-$m$ Uglov data}\}.
\end{align}
In particular, if the latter sends $(K, \vec{r}, \mathbf{a})$ to $(N, \vec{s}, \mathbf{b})$, then $\Upsilon_m^l$ restricts to a bijection 
\begin{align}
	\{|\vec{\pi}, \vec{r}\rangle \mid \vec{\pi} \in \mathbf{a}\} \xrightarrow{\sim} \{|\vec{\varpi}, \vec{s}\rangle \mid \vec{\varpi} \in \mathbf{b}\}.
\end{align}
\end{thm}

Theorem \ref{thm:uglov} is remarkable because there is no direct relationship between the algebras $H_{M, l}(\eta_m, \eta_m^{\vec{r}})$ and $H_{N, m}(\eta_l, \eta_l^{\vec{s}})$ beyond the numerics of their parameters.

\subsection{}\label{subsec:level-rank}

We now sketch how Theorem \ref{thm:uglov} is related to Frenkel's \dfemph{level-rank duality} for affine Lie algebras, or rather, for quantizations $\mathrm{U}_v(\widehat{\mathfrak{sl}}_l)_\bb{Q}$ and $\mathrm{U}_v(\widehat{\mathfrak{sl}}_m)_\bb{Q}$ over the field $\bb{Q}(v)$~\cite{frenkel}.

For any integer $l > 0$ and $m$-charge $\vec{s}$, there is an integrable module over the quantum affine algebra $\mathrm{U}_v(\widehat{\mathfrak{sl}}_m)_\bb{Q}$ called the \dfemph{Fock space of level $l$ and charge $\vec{s}$}.
Its underlying vector space is a formal span of charged $m$-partitions:
\begin{align}
	\Lambda_v^{\vec{s}} = \bigoplus_{\vec{\pi} \in \Pi^m} \bb{Q}(v)|\vec{\pi}, \vec{s}\rangle.
\end{align}
It controls the representation theory of the algebras $H_{N, m}(\eta_l, \eta_l^{\vec{s}})$ through the theorem of Ariki mentioned earlier.
Namely, Ariki identifies the sum of the Grothendieck groups $\mathsf{K}_0(H_{N, m}(\eta_l, \eta_l^{\vec{s}}))$ for $N \geq 0$ with an irreducible, highest-weight submodule of $\Lambda_v^{\vec{s}}$.
If we set $|\vec{s}| = s^{(0)} + \cdots + s^{(m - 1)}$, then for any fixed integer $s$, Uglov's commuting actions arise from the linear isomorphisms
\begin{align}
	\bigoplus_{|\vec{r}| = s} \Lambda_v^{\vec{r}} \xleftarrow{\upsilon_l} \Lambda_v^s \xrightarrow{\upsilon_m^\ast} \bigoplus_{|\vec{s}| = s} \Lambda_v^{\vec{s}}
\end{align}
induced by the maps $\upsilon_l$ and $\upsilon_m^\ast$ such that $\upsilon_m^l = \upsilon_m^\ast \circ \upsilon_l^{-1}$. 
(See the end of \S\ref{sec:partitions}.)

\section{Categories}\label{sec:cats}

\subsection{}

After Uglov, several teams of authors worked to categorify Theorem \ref{thm:uglov}(3) to a statement about highest-weight covers of the module categories of the Hecke algebras $H_{N, m}(\vec{x})$.
To explain this, it is convenient to generalize from the groups $S_{N, m}$ to a finite complex reflection group $C$ with reflection representation $V$.
Generalizing $\Br_{N, m}$, the \dfemph{braid group} of $C$ is the fundamental group $\Br_C \vcentcolon= \pi_1(V^\circ/C)$, where $V^\circ \subseteq V$ is the open locus where $C$ acts freely.
Generalizing $H_{N, m}(\vec{x})$, the \dfemph{Hecke algebra} $H_C(\vec{x})$ is a certain quotient of $\bb{Z}[\vec{x}^{\pm 1}]\Br_C$, where $\vec{x}$ is now a collection of indeterminates indexed in terms of the reflection hyperplanes in $V$ and the orders of certain corresponding complex reflections~\cite{bmr}.
Extending scalars to $\bb{C}$, and specializing $\vec{x}$ to a vector of nonzero complex numbers $\vec{\xi}$, we obtain an algebra $H_C(\vec{\xi})$, a decomposition map $\mathrm{d}_{\vec{\xi}}$ generalizing the map in \eqref{eq:decomposition-map}, and a notion of $\vec{\xi}$-blocks.
The passage from $H_C(\vec{x})$ to $H_C(\vec{\xi})$ may be viewed as fixing conditions on monodromy over $V^\circ/C$.

The \dfemph{rational double affine Hecke algebra} or \dfemph{rational Cherednik algebra} of $C$, at a vector of complex numbers $\vec{\nu}$ with the same indices as $\vec{\xi}$, is an algebra of deformed polynomial \emph{differential operators} over $V \sslash C$.
It takes the form
\begin{align}
D_C^\rat(\vec{\nu}) = (\bb{C}C \ltimes \ur{T}(V\oplus V^\ast))/I(\vec{\nu}),
\end{align}
where $\ur{T}(-)$ denotes ``tensor algebra'' and $I(\vec{\nu})$ is a certain two-sided ideal, symmetrizing $V$ and $V^\ast$, and deforming the Heisenberg--Weyl relations $xy - yx = \langle x, y\rangle$ for $x \in V$ and $y \in V^\ast$.
As shown in~\cite{ggor}, the rational Cherednik algebra shares many features with the universal enveloping algebras of semisimple Lie algebras:
It has a triangular decomposition, where $\bb{C}C$ plays the role of the Cartan subalgebra, and an analogue of the Bernstein--Gelfand--Gelfand category $\mathsf{O}$, which we will denote $\mathsf{O}_C(\vec{\nu})$.
In particular, $\mathsf{O}_C(\vec{\nu})$ is a highest-weight category whose simple objects are indexed by $\Irr(C)$.
For any $\chi \in \Irr(C)$, we write $\Delta_{\vec{\nu}}(\chi)$ to denote the \dfemph{Verma} or \dfemph{standard module} with the corresponding simple quotient.

When $\vec{\xi}$ is related to $\vec{\nu}$ by a certain exponential formula, modules over $H_C(\vec{\xi})$ and $D_C^\rat(\vec{\nu})$ are related through a Riemann--Hilbert correspondence.
Namely, localizing a $D_C^\rat(\vec{\nu})$-module to $V^\circ/C$, then taking monodromy along a Knizhnik--Zamolodchikov-type connection, defines a functor
\begin{align}\label{eq:kz}
	\mathsf{KZ} \colon \mathsf{O}_C(\vec{\nu}) \to \mathsf{Mod}(H_C(\vec{\xi})).
\end{align}
When $C$ is a \emph{real} reflection group, the formula is just $\vec{\xi} = \exp(2\pi i\vec{\nu})$.

Ginzburg--Guay--Opdam--Rouquier show that the $\mathsf{KZ}$ functor induces a bijection between the block \emph{subcategories} of $\mathsf{O}_C(\vec{\nu})$ and those of $\mathsf{Mod}(H_C(\vec{\xi}))$~\cite{ggor}.
In fact, $\sf{KZ}$ is representable by a certain projective object $P$, and the pair $(\mathsf{O}_C(\vec{\nu}), P)$ forms a \dfemph{highest-weight cover} of $\mathsf{Mod}(H_C(\vec{\xi}))$ in the sense of Rouquier~\cite{rouquier_08}.

\subsection{}

Let $(K, \vec{r}, \mathbf{a})$ and $(N, \vec{s}, \mathbf{b})$ be Uglov data corresponding to each other in the setup of \Cref{thm:uglov}.
The following result was essentially conjectured by Chuang--Miyachi~\cite{cm_12}, and proved by Rouquier--Shan--Varagnolo--Vasserot~\cite[Thm.\@ 7.4]{rsvv}. 
(See also~\cite[Cor.\ 6.5]{svv}.)

\begin{thm}[Categorical Level-Rank Duality] 
Let $(K, \vec{r}, \mathbf{a})$ and $(N, \vec{s}, \mathbf{b})$ be Uglov data corresponding to each other in the setup of \Cref{thm:uglov}.
Then the bijection $\mathsf{a} \xrightarrow{\sim} \mathsf{b}$ is categorified by an equivalence of bounded derived categories
	\begin{align}
	\mathsf{D}^b(\mathsf{O}_{S_{K, l}}(\vec{\nu}_l)_\mathbf{a}) \simeq \mathsf{D}^b(\mathsf{O}_{S_{N, m}}(\vec{\nu}_m)_\mathbf{b})
	\end{align}
	that combines Koszul and Ringel dualities.
	Above, the Cherednik parameters $\vec{\nu}_l \in \bb{C}^l$ and $\vec{\nu}_m \in \bb{C}^m$ respectively map to the Hecke parameters $(\eta_m, \eta_m^{\vec{r}})$ and $(\eta_l, \eta_l^{\vec{s}})$ under \eqref{eq:kz}; the subscripts $\mathbf{a}, \mathbf{b}$ indicate the blocks of the Cherednik categories $\mathsf{O}$ lifting the appropriate blocks of the Hecke module categories.
\end{thm}

The proof uses further equivalences between the Cherednik categories and truncated parabolic categories $\mathsf{O}$ of the Lie algebras $\widehat{\mathfrak{sl}}_m$, proved independently by Losev~\cite{losev} and Rouquier--Shan--Varagnolo--Vasserot~\cite{rsvv}.

\section{Finite Reductive Groups}\label{sec:reductive}

\subsection{}

We will propose a generalization of \Cref{thm:uglov}, replacing the groups $S_{N, m}$ with complex reflection groups arising from the representation theory of finite groups of Lie type.
As we will explain, the \emph{general linear} case of our story recovers the case of \Cref{thm:uglov} where $l$ and $m$ are coprime.
We expect that the other classical linear cases also recover \Cref{thm:uglov}---in type $D$, possibly after replacing $S_{N, m}$ by an index-$2$ subgroup.

Fix a prime power $q$ and a (connected, smooth) reductive algebraic group $\bb{G}$ over $\bar{\bb{F}}_q$, equipped with an $\bb{F}_q$-structure corresponding to a Frobenius map $F \colon \bb{G} \to\bb{G}$.
We say that $G = \bb{G}^F$ is a \dfemph{finite reductive group}. 
For instance, if $\bb{G}$ is the general linear group $\bb{GL}_n$ over $\bar{\bb{F}}_q$, then different choices of $F$ produce either $G = \GL_n(\bb{F}_q)$ or $G = \ur{GU}_n(\bb{F}_q)$.

Throughout, we use \textbf{boldface} uppercase letters for spaces over $\bar{\bb{F}}_q$, and ordinary \emph{italics} for their loci of $F$-fixed points.

\subsection{}

By the work of Deligne and Lusztig, the irreducible representations of $G$ can all be obtained from certain induction maps
\begin{align}
	R_\bb{M}^\bb{G} \colon \mathsf{K}_0(M) \to \mathsf{K}_0(G),
	\quad\text{where $\bb{M}$ runs over $F$-stable Levi subgroups of $\bb{G}$},
\end{align}
and $M = \bb{M}^F$.
In more detail, $R_\bb{M}^\bb{G}$ arises from commuting actions of $G$ and $M$ on the compactly-supported \'etale cohomology of algebraic varieties $\bb{Y}_\bb{P}^\bb{G}$ over $\bar{\bb{F}}_q$, now called \dfemph{Deligne--Lusztig varieties}.
The precise definition depends on a choice of parabolic subgroup $\bb{P} \subseteq \bb{G}$ containing $\bb{M}$.
For $\bb{F}_q$ of sufficiently large characteristic, the map $R_\bb{M}^\bb{G}$ is independent of $\bb{P}$.

In fact, to construct the irreducibles of $G$, it suffices to run over maximal tori, not Levis.
An irreducible character is \dfemph{unipotent} if and only if it occurs in $R_\bb{T}^\bb{G}(1)$ for some maximal torus $\bb{T}$, where $1$ is the trivial character.
Following Lusztig, let
\begin{align}
	\Uch(G) = \{\text{unipotent irreducible characters of $G$}\}.
\end{align}
Lusztig shows that $\Uch(G)$ can be indexed in a way independent of $q$, depending only on the Weyl group $W$ of $\bb{G}$ itself~\cite{lusztig_78, lusztig_84}.

\subsection{}

An irreducible character $\mu \in \Irr(M)$ is \dfemph{cuspidal} if and only if it does not occur in $R_\bb{L}^\bb{M}(\lambda)$ for some strictly smaller $F$-stable Levi $\bb{L} \subseteq \bb{M}$ and irreducible $\lambda \in \Irr(L)$.
For a unipotent cuspidal $\mu \in \Uch(M)$, we define the \dfemph{Harish-Chandra series} associated with $(M, \mu)$ to be the set
\begin{align}
	\Uch(G)_{M, \mu} = \{\rho \in \Uch(G) \mid (\rho, R_\bb{M}^\bb{G}(\mu))_G \neq 0\}.
\end{align}
Harish-Chandra observed that the sets $\Uch(G)_{M, \mu}$ are pairwise disjoint and partition $\Uch(G)$ as we run over $\bb{G}$-conjugacy classes of such \dfemph{cuspidal pairs} $(M, \mu)$ for which the underlying Levi $\bb{M} \subseteq \bb{G}$ is $F$-maximally split.
This last condition means that $\bb{M} = Z_\bb{G}(\bb{T})^\circ$ for some $F$-stable torus $\bb{T}$, maximally split over $\bb{F}_q$.

Motivated by ideas from the $\ell$-modular representation theory of $G$ at large primes $\ell$, Brou\'e--Malle--Michel discovered that Harish--Chandra's theory is the $m = 1$ case of a theory that exists for any positive integer $m$~\cite{bmm}.
(To obtain a relation between $\ell$-modular representation theory and what follows, one takes $m$ to be the multiplicative order of $q$ modulo $\ell$.)
To state their results cleanly, they introduce the notion of a \dfemph{generic finite reductive group} that interpolates the groups $G$ as we keep the root datum and its Frobenius automorphism fixed, but vary the prime power $q$.
For simplicity, we will avoid this formalism in our discussion below.

Let $\Phi_m$ denote the $m$th cyclotomic polynomial.
It defines a class of $F$-stable tori $\bb{T} \subseteq \bb{G}$, not necessarily maximal, called \dfemph{$\Phi_m$-tori}:
The corresponding groups $T \subseteq G$ are characterized by the property that $|T|$ is a power of $\Phi_m(q)$.

An $F$-stable Levi subgroup $\bb{M} \subseteq \bb{G}$ is called \dfemph{$\Phi_m$-split} if and only if it takes the form $\bb{M} = Z_\bb{G}(\bb{T})^\circ$ for some $\Phi_m$-torus $\bb{T}$.
For such $\bb{M}$, we say that $\mu \in \Uch(M)$ is \dfemph{$\Phi_m$-cuspidal} if and only if it does not occur in $R_\bb{L}^\bb{M}(\lambda)$ for any smaller $\Phi_m$-split Levi $\bb{L}$ and $\lambda \in \Uch(L)$.
In this case, we say that $(M, \mu)$ is a \dfemph{$\Phi_m$-cuspidal pair}.
By taking $m = 1$, we recover the usual notions of maximally split tori, maximally split Levis, and cuspidal pairs.

Note that $G$ acts by conjugation on the set of $\Phi_m$-cuspidal pairs.
The following result is essentially part (1) of Theorem 3.2 in~\cite{bmm}:

\begin{thm}[Brou\'e--Malle--Michel]
	For any fixed integer $m > 0$, the Harish-Chandra series $\Uch(G)_{M, \mu}$ are disjoint and partition $\Uch(G)$, as $(M, \mu)$ runs over a full set of representatives for the conjugacy classes of $\Phi_m$-cuspidal pairs.
\end{thm}

\subsection{}

In the classical case where $m = 1$, Lusztig observed~\cite{lusztig_77b, lusztig_78} that each Harish-Chandra series $\Uch(G)_{M, \mu}$ has a parametrization of the form
\begin{align}
	\chi_{M, \mu} \colon \Uch(G)_{M, \mu} \xrightarrow{\sim} \Irr(W_{M, \mu}^G),
\end{align}
where $W_{M, \mu}^G$ is the stabilizer of $\mu$ under the action of $W_M^G \vcentcolon= N_G(M)/M$.
We say that $W_M^G$, \emph{resp.}\@ $W_{M, \mu}^G$, is the \dfemph{relative Weyl group} of $M$, \emph{resp.}\@ $(M, \mu)$, in $G$.
These form real reflection groups; for most choices of $\mu$, they coincide.

The bijection $\chi_{M, \mu}$ arises from comparing the $G$-commutant of the cohomology of $Y_\bb{P}^\bb{G}$ to the group algebra of $W_{M, \mu}^G$.
More precisely, Lusztig shows that there is a $\bar{\bb{Q}}[x^{\pm 1/\infty}]$-algebra $H_{M, \mu}^G(x)$, obtained from the algebra $H_{W_{M, \mu}^G}(\vec{x})$ discussed earlier by extending scalars and specializing $\vec{x}$, such that for generic $\xi \in \bb{C}^\times$, we have
\begin{align}\label{eq:hecke-complex}
	\bb{C} \otimes_{\bar{\bb{Q}}[\vec{x}^{\pm 1/\infty}]} H_{M, \mu}^G(\xi) \simeq \bb{C}W_{M, \mu}^G,
	\quad\text{where $x \mapsto \xi$},
\end{align}
but at the same time,
\begin{align}\label{eq:hecke-l-adic}
	\bar{\bb{Q}}_\ell \otimes_{\bar{\bb{Q}}[\vec{x}^{\pm 1/\infty}]} H_{M, \mu}^G(x) \simeq \End_G(\ur{H}_c^\ast(Y_\bb{P}^\bb{G}, \bar{\bb{Q}}_\ell)),
	\quad\text{where $x \mapsto q$}.
\end{align}
Then $\chi_{M, \mu}$ arises from a composition of isomorphisms
\begin{align}
	\mathsf{K}_0(G) \xrightarrow{\sim} \mathsf{K}_0(H_{M, \mu}^G(q)) \xleftarrow{\sim} \mathsf{K}_0(H_{M, \mu}^G(x)) \xrightarrow{\sim} \mathsf{K}_0(W_{M, \mu}^G),
\end{align}
where the first arrow arises from the double-centralizer theorem and the latter two arise from Tits deformation.
In fact, these isomorphisms are isometries with respect to natural multiplicity pairings.

Lusztig observed similar results in the cases where the underlying Levi is not $F$-maximally split, but a maximal torus of Coxeter type~\cite{lusztig_76}.
This condition implies that it is $\Phi_h$-split for some integer $h$, called the (twisted) Coxeter number of $G$.
In these cases, the relative Weyl group is cyclic of order $h$.

Brou\'e--Malle observed that for general $m$, the relative Weyl groups arising from $\Phi_m$-split Levis and $\Phi_m$-cuspidal pairs are always \emph{complex} reflection groups~\cite{bm_93}.
Generalizing Lusztig's results, they define an algebra $H_{M, \mu}^G(x)$ for any $G, m, M, \mu$, using case-by-case formulas, such that \eqref{eq:hecke-complex} still holds, and various numerical properties of $H_{M, \mu}^G(x)$ are consistent with \eqref{eq:hecke-l-adic}.

\begin{conj}[Brou\'e--Malle] 
	\eqref{eq:hecke-l-adic} holds for all $G, m, M, \mu$.
\end{conj}

Although few cases of this conjecture are known beyond the ones that Lusztig established, Brou\'e--Malle--Michel were still able to generalize his parametrizations $\chi_{M, \mu}$ to all $G, m, M, \mu$, by direct and case-by-case constructions.
The following result is part (2) of Theorem 3.2 in~\cite{bmm}.

\begin{thm}[Brou\'e--Malle--Michel]\label{thm:bmm}
For any fixed integer $m > 0$ and $\Phi_m$-cuspidal pair $(M, \mu)$, there is a map
	\begin{align}
		(\varepsilon_{M, \mu}, \chi_{M, \mu}) \colon \Uch(G)_{M, \mu} \to \{\pm 1\} \times \Irr(W_{M, \mu}^G)
	\end{align}
	such that $\varepsilon_{M, \mu}\chi_{M, \mu} \colon \bb{Z}\Uch(G)_{M, \mu} \to \mathsf{K}_0(W_{M, \mu}^G)$ defines an isometry with respect to the natural multiplicity pairings.
	In particular, $\chi_{M, \mu}$ is bijective.
	Moreover:
	\begin{enumerate}
		\item 	These maps are compatible with inclusions of $m$-split Levis.
		They intertwine the maps $R_L^M$ with the ordinary induction maps between relative Weyl groups.
		
		\item   The collection of maps $(\varepsilon_{M, \mu}, \chi_{M, \mu})$ is stable under the conjugation action of $G$ on $\Phi_m$-cuspidal pairs.
		
		
	\end{enumerate}
\end{thm}

\begin{rem}
Brou\'e--Malle--Michel show that if $\rho \in \Uch(G)$ satisfies $\Deg_\rho(\zeta_m) \neq 0$, then $\rho \in \Uch(G)_{T, 1}$, where $T$ represents the conjugacy class of $\Phi_m$-split Levis that are maximal tori, and $1$ denotes the trivial character of $T$.
In this case, $\mathrm{Deg}_\rho(\zeta_m) = \varepsilon_{T, 1}(\rho) \deg \chi_{T, 1}(\rho)$~\cite[76]{broue}.
\end{rem}

\section{Finite General Linear Groups}\label{sec:linear}

\subsection{}

Take $\bb{G} = \bb{GL}_n$ under the standard Frobenius map, so that $G = \GL_n(\bb{F}_q)$.
Then there is the following dictionary between the notions from \S\ref{sec:partitions} and the notions we have just introduced. 

First, recall that every unipotent irreducible character of $\GL_n(\bb{F}_q)$ arises from a principal series attached to an irreducible character of its Weyl group $S_n$.
Restated in the language of Harish-Chandra series:
For any maximally split maximal torus $\bb{A} \subseteq \bb{G}$, we have $\Uch(G)_{A, 1} = \Uch(G)$ and $W_{A, 1}^G \simeq S_n$.
At the same time, $\Irr(S_n)$ is indexed by partitions of size $n$.
Thus, we may regard $\chi_{A, 1}$ as a bijection:
\begin{align}\label{eq:ggl}
	\Uch(G) \simeq \{\rho \in \Pi \mid |\rho| = n\}.
\end{align}
Henceforth, we abuse notation by writing $\rho$ to denote both characters of $G$ and the corresponding partitions.

Next, for any $m \geq 1$, it turns out that every $\Phi_m$-split Levi subgroup of $\bb{G}$ arising from a $\Phi_m$-cuspidal pair takes the form $\bb{M} \simeq {\GL_{n - mN}} \times {\bb{S}_m^N}$, where $\bb{S}_m$ is a torus satisfying $\bb{S}_m(\bb{F}_q) = \bb{F}_{q^m}^\times$.
For such $\bb{M}$, the analogue of \eqref{eq:ggl} restricts to a bijection
\begin{align}
	\{\text{$m$-cuspidals $\mu \in \Uch(M)$}\}
	\simeq \{\mu \in \Pi_{\core{m}} \mid |\mu| = n - mN\}.
\end{align}
Henceforth, we abuse notation by writing $\mu$ to denote both cuspidal characters of $M$ and the corresponding $m$-core partitions.
For a fixed $\mu$, \eqref{eq:ggl} restricts to a bijection
\begin{align}\label{eq:ggl-hc}
	\Uch(G)_{M, \mu}
	&\simeq \{\rho \in \Pi \mid \text{$|\rho| = n$ and $\text{$m$-core}(\rho) = \mu$}\}.
\end{align}
At the same time, it turns out that $W_{M, \mu}^G \simeq W_M^G \simeq S_{N, m}$, giving a bijection
\begin{align}
	\Irr(W_{M, \mu}^G)
	&\simeq \{\vec{\pi} \in \Pi^m \mid |\vec{\pi}| = N\}.
\end{align}
Under these identifications, $\chi_{M, \mu} \colon \Uch(G)_{M, \mu} \xrightarrow{\sim} \Irr(W_{M, \mu}^G)$ is essentially the map sending a partition to its $m$-quotient.
To be precise: $\chi_{M, \mu}(\rho)$ is the ``shifted'' $m$-quotient $\Upsilon_m(\rho, \ell_\mu)$, where $\ell_\mu$ denotes the length of $\mu$ as a partition.

\subsection{}

We can now provide the explicit definition of the Brou\'e--Malle algebra $H_{M, \mu}^G(x)$ for $G = \GL_n(\bb{F}_q)$.
Recall that if $W_{M, \mu}^G \simeq S_{N, m}$ for some $N$, then $H_{M, \mu}^G(x)$ is obtained from $H_{W_{M, \mu}^G}(\vec{x})$ by extending scalars to $\bar{\bb{Q}}[x^{\pm 1/\infty}]$ and specializing $\vec{x}$.

By the core-quotient bijection, the $\bb{Z}^m$-component of $\Upsilon_m(\rho, \ell_\mu)$ is the same vector $\vec{b}_m(\mu)$ for all $\rho \in \Pi^m$.
In the notation of \S\ref{sec:reps}, the map $H_{S_{N, m}}(\vec{x}) \to H_{M, \mu}^G(x)$ sends
\begin{align}\label{eq:gl-params}
\left\{\begin{array}{r@{\:}l}
x_\sigma 
	&\mapsto x^m,\\[0.5ex]
x_{\tau, j}
	&\mapsto \text{$x^{m b_m^{(j)}(\mu) + j}$ for all $j$}.
\end{array}\right.
\end{align}
Note that we always have $b_m^{(0)}(\mu) = 0$.

\section{Conjectures and Results}\label{sec:tx}

\subsection{}

We return to the general setting of \S\ref{sec:reductive}.
We define \dfemph{$H_{M, \mu}^G(\xi)$-blocks} of $\Irr(W_{M, \mu}^G)$ similarly to the $\vec{\xi}$-blocks in \S\ref{sec:reps}, but through a composition of the form
\begin{align}
	\Irr(W_{M, \mu}^G) \to \mathsf{K}_0^+(H_{W_{M, \mu}^G}^\circ(\vec{x})) \xrightarrow{\mathrm{d}_\xi} \mathsf{K}_0^+(H_{M, \mu}^G(\xi)),
\end{align}
where $\mathrm{d}_\xi$ is induced by $H_{W_{M, \mu}^G}(\vec{x}) \to H_{M, \mu}^G(x) \to H_{M, \mu}^G(\xi)$.

Keeping $q, \bb{G}, F$, we now consider a $\Phi_l$-cuspidal pair $(L, \lambda)$ and a $\Phi_m$-cuspidal pair $(M, \mu)$, for separate positive integers $l$ and $m$.
Fix a primitive $l$th root of unity $\zeta_l$ and a primitive $m$th root of unity $\zeta_m$.
Our main conjecture is:

\begin{conj}\label[conj]{conj:level-rank-g}
	Let $\Irr(W_{L, \lambda}^G)_{M, \mu}$ and $\Irr(W_{M, \mu}^G)_{L, \lambda} $ be the images of the maps
	\begin{align}
		\Irr(W_{L, \lambda}^G)
		\xleftarrow{\chi_{L, \lambda}}
		\Uch(G)_{L, \lambda} \cap \Uch(G)_{M, \mu}
		\xrightarrow{\chi_{M, \mu}}
		\Irr(W_{M, \mu}^G).
	\end{align}
Then:
	\begin{enumerate}
		\item
		$\Irr(W_{L, \lambda}^G)_{M, \mu}$ and $\Irr(W_{M, \mu}^G)_{L, \lambda}$ are unions of $H_{L, \lambda}^G(\zeta_m)$- and $H_{M, \mu}^G(\zeta_l)$-blocks, respectively.	
		\item
		The bijection $\Irr(W_{L, \lambda}^G)_{M, \mu} \xrightarrow{\sim} \Irr(W_{M, \mu}^G)_{L, \lambda}$ induced by $\chi_{L, \lambda}$ and $\chi_{M, \mu}$ respects blocks.
		That is, it descends to a bijection $\chi_{M, \mu}^{L, \lambda}$ of the form
		\begin{align}
			\{\text{$H_{L, \lambda}^G(\zeta_m)$-blocks $\mathbf{a} \subseteq \Irr(W_{L, \lambda}^G)_{M, \mu}$}\}
			\xrightarrow{\sim}
			\{\text{$H_{M, \mu}^G(\zeta_l)$-blocks $\mathbf{b} \subseteq \Irr(W_{M, \mu}^G)_{L, \lambda}$}\}.
		\end{align}
		
		\item
		If $\chi_{M, \mu}^{L, \lambda}(\mathbf{a}) = \mathbf{b}$, then the bijection $\mathbf{a} \xrightarrow{\sim} \mathbf{b}$ is categorified by a bounded derived equivalence between appropriate highest-weight covers: that is, an equivalence
		\begin{align}
			\mathsf{D}^b(\mathsf{O}_{W_{L, \lambda}^G}(\vec{\nu}_l)_\mathbf{a})
			\simeq 
			\mathsf{D}^b(\mathsf{O}_{W_{M, \mu}^G}(\vec{\nu}_m)_\mathbf{b})
		\end{align}
		for any vectors $\vec{\nu}_l, \vec{\nu}_m$ related to $\mathbf{a}, \mathbf{b}$ by certain numerical conditions.
		
	\end{enumerate}
\end{conj}

\subsection{}

Our main theorem proves \Cref{conj:level-rank-g} for the groups $G = \GL(\bb{F}_q)$ and coprime $l, m$, via \Cref{thm:uglov}.
To bridge the two settings, we need two further observations.

First, from comparing \eqref{eq:ariki-koike-params} to \eqref{eq:gl-params}, we see that $H_{N, m}(\eta_l, \eta_l^{\vec{s}})$ takes the form $H^G_{M, \mu}(\zeta_l)$ if and only if 
\begin{align}
\eta_l &= \zeta_l^m,\\
m\vec{s} &\equiv m\vec{b}_m(\mu) + (0, 1, \ldots, m - 1) \pmod{l}.
\end{align}
Henceforth, we assume that $l, m$ are coprime, so that $\zeta_l^m$ remains a primitive $l$th root of unity and $\zeta_m^l$ remains a primitive $m$th root of unity.

Second, observe that the affine symmetric group $\bb{Z}^m \rtimes S_m$ acts on the set of $m$-charges.
There is a corresponding action on the set of possible parameters $\vec{\nu}$ for the rational Cherednik algebra, such that if $\vec{\nu}$ maps to $(\eta_l, \eta_l^{\vec{s}}) = (\zeta_l^m, \zeta_l^{m\vec{s}})$ under \eqref{eq:kz}, then $w \cdot \vec{\nu}$ maps to $(\eta_l, \eta_l^{\bar{w} \cdot \vec{s}}) = (\zeta_l^m, \zeta_l^{m(\bar{w} \cdot \vec{s})})$ under \eqref{eq:kz}, for any $w \in \bb{Z}^m \rtimes S_m$ with finite part $\bar{w} \in S_m$.
Rouquier conjectured, and Gordon--Losev~\cite[\S{5.13}--{5.14}]{gl} and Webster~\cite[Thm.\@ B(3)]{webster} independently proved, equivalences of the form $\mathsf{D}^b(\mathsf{O}_{S_{N, m}}(\vec{\nu})) \xrightarrow{\sim} \mathsf{D}^b(\mathsf{O}_{S_{N, m}}(w \cdot \vec{\nu}))$.
As we run over $w$, these functors define a strong categorical action of the affine braid group on the set of such categories.

The proof of the following theorem relies on the abacus combinatorics discussed in \S\ref{sec:partitions}, as well as the classification of $\vec{\xi}$-blocks due to Lyle--Mathas.
In \cite{tx}, see Theorem 3, Theorem 4, Corollary 5, and their proofs in \S{6--7}.

\begin{thm}\label[thm]{thm:main}
Suppose that $G = \GL_n(\bb{F}_q)$ and that $l$ and $m$ are coprime.
	Then all three parts of \Cref{conj:level-rank-g} hold.
	Moreover:
	\begin{enumerate}
		\item 	$\Irr(W_{L, \lambda}^G)_{M, \mu}$ and $\Irr(W_{M, \mu}^G)_{L, \lambda}$ are single $H^G_{L, \lambda}(\zeta_m)$- and $H^G_{M, \mu}(\zeta_l)$-blocks, respectively.
		
		\item 	The bijections in part (2) of \Cref{conj:level-rank-g} are essentially the bijections $\Upsilon_m^l$, in that there are affine permutations $w_l \in S_l \ltimes \bb{Z}^l$ and $w_m \in S_m \ltimes \bb{Z}^m$ giving a commutative diagram:
		\begin{equation}
			\begin{tikzpicture}[baseline=(current bounding box.center), >=stealth]
				\matrix(m)[matrix of math nodes, row sep=2.5em, column sep=2.5em, text height=2ex, text depth=0.5ex]
				{ 		
					\Irr(W_{L, \lambda}^G)
					&\Uch(G)_{L, \lambda} \cap \Uch(G)_{M, \mu}
					&\Irr(W_{M, \mu}^G)\\
					\Pi^l \times \bb{Z}^l
					&\Pi
					&\Pi^m \times \bb{Z}^m\\	
					\Pi^l \times \bb{Z}^l
					&
					&\Pi^m \times \bb{Z}^m\\
				};
				\path[->, font=\scriptsize, auto]
				(m-1-2)
				edge node[above]{$\chi_{L, \lambda}$} (m-1-1)
				edge node{$\chi_{M, \mu}$} (m-1-3)
				edge (m-2-2)
				(m-1-1)
				edge (m-2-1)
				(m-1-3)
				edge (m-2-3)
				(m-2-2)
				edge node[above]{$\Upsilon_l(\rho, \ell_\lambda) \mapsfrom \rho$} (m-2-1)
				edge node{$\rho \mapsto \Upsilon_m(\rho, \ell_\mu)$} (m-2-3)
				(m-2-1)
				edge node[left]{$w_l$} (m-3-1)
				(m-2-3)
				edge node{$w_m$} (m-3-3)
				(m-3-1)
				edge node{$\Upsilon_m^l$} (m-3-3);
			\end{tikzpicture}
		\end{equation}
		Above, $\ell_\lambda, \ell_\mu$ refer to the lengths of $\lambda$ and $\mu$ as partitions.
		
		\item 	
			We can choose $w_m$ so that
			\begin{align} 
			m(w_m \cdot \vec{b}_m(\mu)) \equiv m\vec{b}_m(\mu) + (0, 1, \ldots, m - 1) \pmod{l},
			\end{align} 
		and similarly with $l, m, \lambda, w_l$ in place of $m, l, \mu, w_m$.
		
	\end{enumerate}
\end{thm}

We expect that, up to modifiying both parts of \Cref{thm:main}, we can remove the coprimality hypothesis from the corollary.
Below, \Cref{ex:multiple-blocks} shows that if $G = \GL_n(\bb{F}_q)$, but $l$ and $m$ are not coprime, then multiple blocks can appear in both $\Irr(W_{L, \lambda}^G)_{M, \mu}$ and $\Irr(W_{M, \mu}^G)_{L, \lambda}$.

We checked in \cite[\S{8}]{tx} that in almost all cases where $\bb{G}$ is exceptional, \Cref{conj:level-rank-g}(1) is compatible with the sizes of the sets $\Uch(G)_{L, \lambda} \cap \Uch(G)_{M, \mu}$ and of the blocks appearing in $\Irr(W_{L, \lambda}^G)_{M, \mu}$ and $\Irr(W_{M, \mu}^G)_{L, \lambda}$, where these can be deduced.
These cases provide further examples where $\Irr(W_{L, \lambda}^G)_{M, \mu}$ and $\Irr(W_{M, \mu}^G)_{L, \lambda}$ each consist of multiple blocks.
Later, after \cite{tx} appeared on arXiv, Chlouveraki--Malle proved \Cref{conj:level-rank-g}(1) in almost all of these exceptional-type cases~\cite{cm_zeitschrift}.

\begin{ex}
	Let $G = \GL_8(\bb{F}_q)$ and $l = 4$ and $m = 3$.
	Let $\lambda = (2^2)$, a $4$-core, and $\mu = (2)$, a $3$-core.
	Then \eqref{eq:ggl-hc} gives
	\begin{align}
		\Uch(G)_{L, \lambda} 
		&\simeq \{(6, 2), (5, 3), (2^3, 1^2), (2^2, 1^4)\},\\
		\Uch(G)_{M, \mu}
		&\simeq \{(8), (5, 3), (5, 2, 1), (5, 1^3), (4, 3, 1), (3^2, 1^2), (2^4), (2^2, 1^4), (2, 1^6)\},
	\end{align} 
	from which $\Uch(G)_{L, \lambda} \cap \Uch(G)_{M, \mu}
	\simeq \{(5, 3), (2^2, 1^4)\}$.
	
	Let $\pi = (5, 3)$.
	The abaci that represent $\Upsilon_l(\pi, \ell_\lambda) = \Upsilon_4((5, 3), 2)$ and $\Upsilon_m(\pi, \ell_\mu) = \Upsilon_3((5, 3), 1)$ are below.
	\begin{align}
		\small
		\begin{array}{@{\quad\cdots\quad}cccccc}
			-12 &-8 &-4 &\erase{-}- &\erase{-}- \\
			-11 &-7 &-3 &\erase{-}- &\erase{-}- \\
			-10 &-6 &-2 &\erase{-}- &\erase{-}6 \\
			-9 &-5 &-1 &\erase{-}3 &
		\end{array}
		\qquad
		\begin{array}{@{\quad\cdots\quad}cccccccc}
			-12 &-9 &-6 &-3 &\erase{-}- &\erase{-}- & &\\
			-11 &-8 &-5 &-2 &\erase{-}- &\erase{-}- & &\\
			-10 &-7 &-4 &\erase{-}- &\erase{-}2 &\erase{-}5 & &
		\end{array}
	\end{align}
	Let $w_4 = (1, 0, 0, -1) \circ {\bar{w}_4}$, where $\bar{w}_4 = [1032] \in S_4$.
	Let $w_3 = (2, 0, -1) \circ {\bar{w}_3}$, where $\bar{w}_3 = [201] \in S_3$.
	Then $w_3$ satisfies the congruence in \Cref{thm:main}(3), and similarly with $w_4$.
	Under these affine permutations, the abaci above become the abaci in \Crefrange{ex:3-abacus}{ex:4-abacus}.
	\begin{align}
		\small
		\begin{array}{@{\quad\cdots\quad}cccccc}
			-12 &-8 &-4 &\erase{-}0 \\
			-11 &-7 &-3 &\erase{-}- \\
			-10 &-6 &-2 &\erase{-}2 \\
			-9 &-5 &\erase{-}- &\erase{-}3 
		\end{array}
		\qquad
		\begin{array}{@{\quad\cdots\quad}cccccccc}
			-12 &-9 &-6 &-3 &\erase{-}0 &\erase{-}- &\erase{-}6 &\erase{-}9\\
			-11 &-8 &-5 &-2 &\erase{-}- &\erase{-}- &\erase{-}- &\\
			-10 &-7 &-4 &\erase{-}- &\erase{-}- &\erase{-}- &\erase{-}- &
		\end{array}
	\end{align}
	We saw in \Cref{ex:4-abacus} that $\Upsilon_4^3$ transforms the left abacus into the right one.
\end{ex}

\begin{ex}\label{ex:multiple-blocks}
We demonstrate that if $G = \GL_n(\bb{F}_q)$, but $\ell, m$ are not coprime, then multiple blocks can appear in $\Irr(W_{L, \lambda}^G)_{M, \mu}$ and $\Irr(W_{M, \mu}^G)_{L, \lambda}$. 

Let $n = 4$ and $l = 2$ and $m = 4$.
Let $(L, \lambda)$ be the $\Phi_2$-cuspidal pair corresponding to the empty $2$-core.
Let $(M, \mu)$ be the $\Phi_4$-cupsidal pair corresponding to the $4$-core $(1^2)$.
Then $W^G_{L, \lambda} \simeq S_{3, 2}$ and $W^G_{L, \lambda} \simeq S_{1,4}$.
Under the bijection~\eqref{eq:ggl-hc},
\begin{align}
\Uch(G)_{M, \mu} \cap \Uch(G)_{L,\lambda}
	\simeq \{(1^6), (4, 2), (5, 1), (2^3)\}.
\end{align}
The maps $\chi_{L, \lambda}, \chi_{M, \mu}$ are given by:
\begin{align}
\begin{array}{l|ll}
\rho &\chi_{L, \lambda}(\rho) &\chi_{M, \mu}(\rho)\\
\hline 
(1^6)
	&(\emptyset, (1^3))
	&((1), \emptyset, \emptyset, \emptyset)\\
(4, 2)
	&((2), (1))
	&(\emptyset, (1), \emptyset, \emptyset)\\
(5, 1)
	&(\emptyset, (3))
	&(\emptyset, \emptyset, (1), \emptyset)\\
(2^3)
	&((1^2), (1))
	&(\emptyset, \emptyset, \emptyset, (1))
\end{array}
\end{align}
By~\eqref{eq:gl-params}, $H^G_{L, \lambda}(x)$ has the Hecke relations $(\sigma_i - 1)(\sigma_i + x^2) = (\tau - 1)(\tau - x^3) = 0$, while $H^G_{M, \mu}(x)$ has the single relation $(\tau - 1)(\tau - x^3)(\tau - x^6)(\tau - x^3) = 0$.
One can check $\Irr(W_{L, \lambda}^G)_{M, \mu}$ decomposes into the two $H^G_{L,\lambda}(\zeta_4)$-blocks
\begin{align} 
\{(\emptyset, (1^3)), (\emptyset, (3))\},\quad 
\{((2), (1)), ((1^2), (1))\},
\end{align}
while $\Irr(W_{M, \mu}^G)_{L, \lambda}$ decomposes into the two $H^G_{M, \mu}(\zeta_2)$-blocks
\begin{align} 
\{((1), \emptyset, \emptyset, \emptyset), (\emptyset, \emptyset, (1), \emptyset)\},\quad
\{(\emptyset, (1), \emptyset, \emptyset), 
(\emptyset, \emptyset, \emptyset, (1))\}.
\end{align} 
One sees that these decompositions satisfy part (2) of \Cref{conj:level-rank-g}.
\end{ex}

\section{Postscript}

An important theme that our discussion has omitted:
When $\bb{L}$ and $\bb{M}$ are maximal tori, we expect novel \emph{geometric} incarnations of the duality in \Cref{conj:level-rank-g}, via bimodules constructed from the cohomology of algebraic varieties.
They suggest analogies between these varieties and certain $F$-twisted Steinberg varieties related to Deligne--Lusztig theory, but with roots of unity replacing the prime power $q$.

One family of incarnations is of Dolbeault nature:
The (ind-)varieties are \dfemph{homogeneous affine Springer fibers} that first appeared in work of Lusztig--Smelt, and more recently, in work of Varagnolo--Vasserot~\cite{vv}, Oblomkov--Yun~\cite{oy}, and Boixeda-Alvarez--Losev--Kivinen~\cite{bl} on the geometric representation theory of (degenerate) double affine Hecke algebras.
Another family is of Betti nature:
The varieties are \dfemph{braid Steinberg varieties} that first appeared in work of Trinh on triply-graded link homology~\cite{trinh_21}.


\bibliographystyle{alphaurl}
\bibliography{blocks}

@article {bl,
    AUTHOR = {Boixeda Alvarez, P. and Losev, I.},
     TITLE = {Affine {S}pringer fibers, {P}rocesi bundles, and {C}herednik  algebras},
      NOTE = {With an appendix by Boixeda Alvarez, Losev and O. Kivinen},
   JOURNAL = {Duke Math. J.},
  FJOURNAL = {Duke Mathematical Journal},
    VOLUME = {173},
      YEAR = {2024},
    NUMBER = {5},
     PAGES = {807--872},
   MRCLASS = {16G99},
  MRNUMBER = {4739836},
       DOI = {10.1215/00127094-2023-0027},
}

@incollection {broue,
    AUTHOR = {Brou\'{e}, M.},
     TITLE = {Reflection groups, braid groups, {H}ecke algebras, finite reductive groups},
 BOOKTITLE = {Current developments in mathematics, 2000},
     PAGES = {1--107},
 PUBLISHER = {Int. Press, Somerville, MA},
      YEAR = {2001},
   MRCLASS = {20F55 (20C08 20F36 20G05)},
  MRNUMBER = {1882533},
}

@incollection {bm_93,
    AUTHOR = {Brou\'e, Michel and Malle, Gunter},
     TITLE = {Zyklotomische {H}eckealgebren},
 BOOKTITLE = {Repr\'esentations unipotentes g\'en\'eriques et blocs des groupes r\'eductifs finis},
   JOURNAL = {Ast\'erisque},
  FJOURNAL = {Ast\'erisque},
    NUMBER = {212},
      YEAR = {1993},
     PAGES = {119--189},
      ISSN = {0303-1179,2492-5926},
   MRCLASS = {20G40 (20C20 20H15)},
  MRNUMBER = {1235834},
MRREVIEWER = {Fran\c cois\ Digne},
}

@incollection {bmm,
    AUTHOR = {Brou\'{e}, M. and Malle, G. and Michel, J.},
     TITLE = {Generic blocks of finite reductive groups},
 BOOKTITLE = {Repr\'esentations unipotentes g\'en\'eriques et blocs des groupes r\'eductifs finis},
   JOURNAL = {Ast\'{e}risque},
  FJOURNAL = {Ast\'{e}risque},
    NUMBER = {212},
      YEAR = {1993},
     PAGES = {7--92},
      ISSN = {0303-1179},
   MRCLASS = {20G05 (20C20)},
  MRNUMBER = {1235832},
MRREVIEWER = {Bhama Srinivasan},
}

@article {bmr,
    AUTHOR = {Brou\'{e}, M. and Malle, G. and Rouquier, R.},
     TITLE = {Complex reflection groups, braid groups, {H}ecke algebras},
   JOURNAL = {J. Reine Angew. Math.},
  FJOURNAL = {Journal f\"{u}r die Reine und Angewandte Mathematik. [Crelle's Journal]},
    VOLUME = {500},
      YEAR = {1998},
     PAGES = {127--190},
      ISSN = {0075-4102},
   MRCLASS = {20F55 (20C05 20F36)},
  MRNUMBER = {1637497},
MRREVIEWER = {Arjeh M. Cohen},
}

@article{cm_zeitschrift,
	author = {Chlouveraki, Maria and Malle, Gunter},
	title = {Intersections of blocks of cyclotomic {Hecke} algebras},
	fjournal = {Mathematische Zeitschrift},
	journal = {Math. Z.},
	issn = {0025-5874},
	volume = {312},
	number = {4},
	pages = {17},
	note = {Id/No 99},
	year = {2026},
	doi = {10.1007/s00209-026-03986-w},
}

@misc {cm_12,
    AUTHOR = {Chuang, J. and Miyachi, H.},
     TITLE = {Hidden {H}ecke Algebras and {K}oszul Duality},
      YEAR = {2012},
       URL = {https://www.math.nagoya-u.ac.jp/~miyachi/preprints/lrkszl21.pdf}
}

@incollection {frenkel,
    AUTHOR = {Frenkel, I. B.},
     TITLE = {Representations of affine {L}ie algebras, {H}ecke modular forms and {K}orteweg--de {V}ries type equations},
 BOOKTITLE = {Lie algebras and related topics ({N}ew {B}runswick, {N}.{J}., 1981)},
    SERIES = {Lecture Notes in Math.},
    VOLUME = {933},
     PAGES = {71--110},
 PUBLISHER = {Springer, Berlin--New York},
      YEAR = {1982},
   MRCLASS = {17B35 (10D40 35Q20 58F07 81G20)},
  MRNUMBER = {675108},
MRREVIEWER = {Jouko Mickelsson},
}

@book {gj,
    AUTHOR = {Geck, M. and Jacon, N.},
     TITLE = {Representations of {H}ecke algebras at roots of unity},
    SERIES = {Algebra and Applications},
    VOLUME = {15},
 PUBLISHER = {Springer-Verlag London, Ltd., London},
      YEAR = {2011},
     PAGES = {xii+401},
      ISBN = {978-0-85729-715-0},
   MRCLASS = {20C08 (17B37)},
  MRNUMBER = {2799052},
MRREVIEWER = {Maria Chlouveraki},
       DOI = {10.1007/978-0-85729-716-7},
}

@book {gp,
    AUTHOR = {Geck, M. and Pfeiffer, G.},
     TITLE = {Characters of finite {C}oxeter groups and {I}wahori--{H}ecke algebras},
    SERIES = {London Mathematical Society Monographs. New Series},
    VOLUME = {21},
 PUBLISHER = {The Clarendon Press, Oxford University Press, New York},
      YEAR = {2000},
      ISBN = {0-19-850250-8},
   MRCLASS = {20C15 (20C08 20F55)},
  MRNUMBER = {1778802},
MRREVIEWER = {Jian-yi Shi},
}

@article {ggor,
    AUTHOR = {Ginzburg, V. and Guay, N. and Opdam, E. and Rouquier, R.},
     TITLE = {On the category {$\mathcal{O}$} for rational {C}herednik algebras},
   JOURNAL = {Invent. Math.},
  FJOURNAL = {Inventiones Mathematicae},
    VOLUME = {154},
      YEAR = {2003},
    NUMBER = {3},
     PAGES = {617--651},
   MRCLASS = {20C08 (20F55)},
  MRNUMBER = {2018786},
MRREVIEWER = {Iain G. Gordon},
       DOI = {10.1007/s00222-003-0313-8},
}

@article {gl,
    AUTHOR = {Gordon, Iain G. and Losev, Ivan},
     TITLE = {On category {$\mathcal{O}$} for cyclotomic rational {C}herednik algebras},
   JOURNAL = {J. Eur. Math. Soc. (JEMS)},
  FJOURNAL = {Journal of the European Mathematical Society (JEMS)},
    VOLUME = {16},
      YEAR = {2014},
    NUMBER = {5},
     PAGES = {1017--1079},
      ISSN = {1435-9855,1435-9863},
   MRCLASS = {17B10 (16E35 20C08 53D55)},
  MRNUMBER = {3210960},
MRREVIEWER = {Selvaraj\ Chelliah},
       DOI = {10.4171/JEMS/454},
}

@book {jk,
    AUTHOR = {James, G. and Kerber, A.},
     TITLE = {The representation theory of the symmetric group},
    SERIES = {Encyclopedia of Mathematics and its Applications},
    VOLUME = {16},
      NOTE = {With a foreword by P. M. Cohn, With an introduction by Gilbert de B. Robinson},
 PUBLISHER = {Addison-Wesley Publishing Co., Reading, MA},
      YEAR = {1981},
     PAGES = {xxviii+510},
      ISBN = {0-201-13515-9},
   MRCLASS = {20-02 (20C30)},
  MRNUMBER = {644144},
MRREVIEWER = {A. O. Morris},
}

@article {losev,
    AUTHOR = {Losev, I.},
     TITLE = {Proof of {V}aragnolo--{V}asserot conjecture on cyclotomic categories {$\mathcal{O}$}},
   JOURNAL = {Selecta Math. (N.S.)},
  FJOURNAL = {Selecta Mathematica. New Series},
    VOLUME = {22},
      YEAR = {2016},
    NUMBER = {2},
     PAGES = {631--668},
   MRCLASS = {20C08 (16G99 17B10 17B67)},
  MRNUMBER = {3477332},
MRREVIEWER = {Volodymyr Mazorchuk},
       DOI = {10.1007/s00029-015-0209-7},
}

@article {lusztig_76,
    AUTHOR = {Lusztig, G.},
     TITLE = {On the finiteness of the number of unipotent classes},
   JOURNAL = {Invent. Math.},
  FJOURNAL = {Inventiones Mathematicae},
    VOLUME = {34},
      YEAR = {1976},
    NUMBER = {3},
     PAGES = {201--213},
   MRCLASS = {20G40},
  MRNUMBER = {419635},
MRREVIEWER = {R. Steinberg},
       DOI = {10.1007/BF01403067},
}

@article {lusztig_77b,
	AUTHOR = {Lusztig, G.},
	TITLE = {Irreducible representations of finite classical groups},
	JOURNAL = {Invent. Math.},
	FJOURNAL = {Inventiones Mathematicae},
	VOLUME = {43},
	YEAR = {1977},
	NUMBER = {2},
	PAGES = {125--175},
	MRCLASS = {20C15},
	MRNUMBER = {463275},
	MRREVIEWER = {James E. Humphreys},
	DOI = {10.1007/BF01390002},
}

@book {lusztig_78,
    AUTHOR = {Lusztig, G.},
     TITLE = {Representations of finite {C}hevalley groups},
    SERIES = {CBMS Regional Conference Series in Mathematics},
    VOLUME = {39},
      NOTE = {Expository lectures from the CBMS Regional Conference held at Madison, Wis., August 8--12, 1977},
 PUBLISHER = {American Mathematical Society, Providence, RI},
      YEAR = {1978},
      ISBN = {0-8218-1689-6},
   MRCLASS = {20G05},
  MRNUMBER = {518617},
MRREVIEWER = {David B. Surowski},
}

@book {lusztig_84,
    AUTHOR = {Lusztig, G.},
     TITLE = {Characters of reductive groups over a finite field},
    SERIES = {Annals of Mathematics Studies},
    VOLUME = {107},
 PUBLISHER = {Princeton University Press, Princeton, NJ},
      YEAR = {1984},
     PAGES = {xxi+384},
      ISBN = {0-691-08350-9; 0-691-08351-7},
   MRCLASS = {20G05 (14L20 20C15)},
  MRNUMBER = {742472},
MRREVIEWER = {Bhama Srinivasan},
       DOI = {10.1515/9781400881772},
}

@article {oy,
    AUTHOR = {Oblomkov, A. and Yun, Z.},
     TITLE = {Geometric representations of graded and rational {C}herednik algebras},
   JOURNAL = {Adv. Math.},
  FJOURNAL = {Advances in Mathematics},
    VOLUME = {292},
      YEAR = {2016},
     PAGES = {601--706}, 
   MRCLASS = {20C08 (14L35 20G25)},
  MRNUMBER = {3464031},
MRREVIEWER = {Simon Riche},
       DOI = {10.1016/j.aim.2016.01.015},
}

@article {rouquier_08,
    AUTHOR = {Rouquier, R.},
     TITLE = {{$q$}-{S}chur algebras and complex reflection groups},
   JOURNAL = {Mosc. Math. J.},
  FJOURNAL = {Moscow Mathematical Journal},
    VOLUME = {8},
      YEAR = {2008},
    NUMBER = {1},
     PAGES = {119--158, 184},
   MRCLASS = {20G43 (16G30 20C08 20F55)},
  MRNUMBER = {2422270},
MRREVIEWER = {Andrew Mathas},
       DOI = {10.17323/1609-4514-2008-8-1-119-158},
}

@article {rsvv,
    AUTHOR = {Rouquier, R. and Shan, P. and Varagnolo, M. and Vasserot, E.},
     TITLE = {Categorifications and cyclotomic rational double affine {H}ecke algebras},
   JOURNAL = {Invent. Math.},
  FJOURNAL = {Inventiones Mathematicae},
    VOLUME = {204},
      YEAR = {2016},
    NUMBER = {3},
     PAGES = {671--786},
   MRCLASS = {20C08 (18D10 20G43)},
  MRNUMBER = {3502064},
MRREVIEWER = {Nicolas Jacon},
       DOI = {10.1007/s00222-015-0623-7},
}

@article {svv,
    AUTHOR = {Shan, P. and Varagnolo, M. and Vasserot, E.},
     TITLE = {{K}oszul duality of affine {K}ac--{M}oody algebras and cyclotomic rational double affine {H}ecke algebras},
   JOURNAL = {Adv. Math.},
  FJOURNAL = {Advances in Mathematics},
    VOLUME = {262},
      YEAR = {2014},
     PAGES = {370--435},
   MRCLASS = {17B67},
  MRNUMBER = {3228432},
MRREVIEWER = {Weiqiang Wang},
       DOI = {10.1016/j.aim.2014.05.012},
}

@misc {trinh_21,
    AUTHOR = {Trinh, M. Q.},
     TITLE = {From the {H}ecke Category to the Unipotent Locus}, 
      YEAR = {2021},
    EPRINT = {2106.07444},
EPRINTTYPE = {arXiv},
      NOTE = {v2},
}

@misc {tx,
	AUTHOR = {Trinh, M. Q. and Xue, Ting},
	TITLE = {Level-Rank Dualities from {$\Phi$}-Cuspidal Pairs and Affine {S}pringer Fibers}, 
	YEAR = {2023},
	EPRINT = {2311.17106},
	EPRINTTYPE = {arXiv},
	NOTE = {v2},
}

@incollection {uglov,
    AUTHOR = {Uglov, D.},
     TITLE = {Canonical bases of higher-level {$q$}-deformed {F}ock spaces and {K}azhdan--{L}usztig polynomials},
 BOOKTITLE = {Physical combinatorics ({K}yoto, 1999)},
    SERIES = {Progr. Math.},
    VOLUME = {191},
     PAGES = {249--299},
 PUBLISHER = {Birkh\"{a}user Boston, Boston, MA},
      YEAR = {2000},
   MRCLASS = {17B37 (17B67)},
  MRNUMBER = {1768086},
MRREVIEWER = {Nicol\'{a}s Andruskiewitsch},
}

@article {vv,
    AUTHOR = {Varagnolo, M. and Vasserot, E.},
     TITLE = {Finite-dimensional representations of {DAHA} and affine {S}pringer fibers: the spherical case},
   JOURNAL = {Duke Math. J.},
  FJOURNAL = {Duke Mathematical Journal},
    VOLUME = {147},
      YEAR = {2009},
    NUMBER = {3},
     PAGES = {439--540},
   MRCLASS = {20C08 (14M15)},
  MRNUMBER = {2510742},
MRREVIEWER = {Maarten Sander Solleveld},
       DOI = {10.1215/00127094-2009-016},       
}

@article {webster,
    AUTHOR = {Webster, B.},
     TITLE = {Rouquier's conjecture and diagrammatic algebra},
   JOURNAL = {Forum Math. Sigma},
  FJOURNAL = {Forum of Mathematics. Sigma},
    VOLUME = {5},
      YEAR = {2017},
     PAGES = {Paper No. e27, 71},
   MRCLASS = {17B37 (16G99 16S37 17B10 17B67 18D99)},
  MRNUMBER = {3732238},
MRREVIEWER = {Jun Hu},
       DOI = {10.1017/fms.2017.17},  
}

\end{document}